\documentclass[12pt,a4paper]{article}
\usepackage{amsmath, amssymb, theorem, latexsym}
\newtheorem{theorem}{Theorem}[section]
\newtheorem{lemma}[theorem]{Lemma}
\newtheorem{proposition}[theorem]{Proposition}
\newtheorem{definition}[theorem]{Definition}

\newtheorem{corollary}[theorem]{Corollary}

\newtheorem{remark}[theorem]{Remark}

\allowdisplaybreaks
\renewcommand{\thefootnote}{\fnsymbol{footnote}}

\newcommand\finbox{~\hfill$\Box$}%

\def\Om {{\Omega}}
\def\la {{\lambda}}

\def \rz { {\mathbb R}}

\def\Og {{\cal O}} 

\def \rz {{\mathbb R}}

\def \Inte{{\rm Int\,}}

\newcommand {\pa}{\partial}

\numberwithin{equation}{section}
\begin{document}
{\centering
\bfseries
{\Large
Nodal Domain Theorems \`a la Courant.}
\\[2\baselineskip]%
{\renewcommand\thefootnote{}%
\footnote{1991 Mathematics Subject Classification 35B05}}
\par
\mdseries
\scshape
\small
A. Ancona$\,^1$,\\
B. Helffer$\,^1$,\\
and T. Hoffmann-Ostenhof$\,^{2,\;3,}$
\footnote{Supported by Ministerium f\"ur Bildung, Wissenschaft und
Kunst der Republik \"Osterreich}~ \\
\par
\upshape
D\a'epartement de Math\a'ematiques, Universit\a'e Paris-Sud$^1$\\
Institut f\"ur Theoretische Chemie, Universit\"at Wien$^2$\\
International Erwin Schr\"odinger Institute for Mathematical Physics$^3$\\

\today

}
\begin{abstract}
Let $H(\Om_0)=-\Delta+V$ be a Schr\"odinger operator on a bounded
domain $\Om_0\subset \mathbb R^d$
with Dirichlet boundary conditions. Suppose  that the $\Om_\ell$ ($\ell \in \{1,\dots,k\}$) are some
pairwise disjoint subsets of $\Om_0$ and that $H(\Om_\ell)$ are the
corresponding Schr\"odinger operators again with Dirichlet boundary
conditions. We investigate the relations between the spectrum of $H(\Om_0)$ and
   the spectra of the $H(\Om_\ell)$. In particular, we derive some inequalities
for the associated spectral counting functions which can be interpreted  as
generalizations of Courant's nodal Theorem. For the case that equality is
achieved we prove converse results. In particular, we use potential
theoretic methods to relate the  $\Om_\ell$  to the nodal domains of
some
eigenfunction of $H(\Omega_0)$.
   \end{abstract}
\section{Introduction}

Consider a Schr\"odinger operator
\begin{equation}\label{HSchr}
    H=-\Delta +V
\end{equation}
on a bounded domain
$\Omega_0 \subset \mathbb R^d$ with Dirichlet boundary conditions.
Further we assume  that $V$ is real valued and satisfies
   $V\in L^\infty(\Om_0)$.
(We could relax this condition and extend our results to
the case  $V\in L^\beta(\Omega_0) $ for some $\beta>d/2$ \cite{JeKe})

$H$ is   selfadjoint if viewed  as the Friedrichs extension
of the quadratic forms of $H$  with form domain
$W_0^{1,2}(\Om_0)$ and form core $C_0^\infty(\Om_0)$ and we
denote it by
$H(\Omega_0)$.
Further  $H(\Omega_0)$  has  compact resolvent. So the spectrum of
$H(\Omega_0)$, $\sigma\big(H(\Om_0)\big)$, can be described by an increasing
sequence of eigenvalues
\begin{equation}\label{laj}
    \lambda_1<\lambda_2\le\lambda_3\le\dots\le \dots
\end{equation}
tending to $+\infty$, such that the associated eigenfunctions $u_k$ form
an orthonormal basis of $L^2(\Om_0)$. $\la_1$ is simple and the corresponding
eigenfunction $u_1$ can be chosen to satisfy, see e.g. \cite{ReedSimon:1978},
\begin{equation}\label{u>0}
\; u_1>0\text{ for all } x\in \Om_0.
\end{equation}

We can assume that the eigenfunctions $u_k$ are real valued and by
elliptic regularity, \cite{GT:1983}(Corollary 8.36),
$u_k$ belongs to $ C^{1,\alpha}(\Om_0)$ for every  $\alpha<1$.
We shall often call for a bounded domain $D$,
$H(D)$,  the corresponding
selfadjoined operator with Dirichlet boundary conditions on $\pa D$.
Its lowest eigenvalue will be denoted by $\la(D)$.

We  denote the zeroset of an eigenvector $u$ by
\begin{equation}\label{Nuk}
N(u)=\overline{\{x\in\Omega_0\,|\; u(x)=0\}}.
\end{equation}
The nodal domains of  $u$, which are by definition the connected
components of $\Omega_0 \setminus N(u)$,
    will be denoted by $D_j, j=1,\dots, \mu(u)$,
where $\mu(u)$ denotes the number of nodal domains of $u$.

Suppose that the $\Om_\ell$ ($\ell=1,2,\dots,k$) are $k$ open pairwise
disjoint subsets
of $\Om_0$. In this paper we shall investigate relations between the
spectrum of $H(\Om_0)$ and the spectra of the $H(\Om_\ell)$. Roughly speaking,
we shall derive an inequality between the counting function of $H(\Om_0)$ and
those of the $H(\Om_\ell)$. This  inequality can be interpreted
as a generalization of Courant's classical nodal domain theorem.
For the case that equality is achieved this will lead to a partial
characterization of the $\Om_\ell$ which will turn out to be related to
the nodal domains of one of the eigenfunction of $H(\Om_0)$.

These results will be given in sections \ref{Section2} and 
\ref{Section3}. From these results some
natural questions of potential theoretic nature arise which will be
analyzed and  answered in section \ref{ancona}.

The proofs of the results stated in sections \ref{Section2} and 
\ref{Section3} are given in sections \ref{Section4}
and \ref{Section5}. In section \ref{Section6} some illustrative 
explicit examples are given.

\section{Main results}\label{Section2}
We start with a result which will turn out to be a  generalization of
Courant's nodal Theorem.
We consider again \eqref{HSchr} on a bounded domain $\Om_0$ and
the corresponding
eigenfunctions and eigenvalues. We first introduce
\begin{equation}\label{nover}
\overline{n}(\lambda,\Omega_0)= \# \{  j\;|\;
\lambda_j(\Omega_0) \leq \lambda\}\;,
\end{equation}
where $\lambda_j(\Omega_0)$ is the $j$-th eigenvalue of $H(\Omega_0)$.\\
We also define
\begin{equation}\label{nunder}
\underline{n} (\lambda, \Omega_0)=\# \{ j\;|\; \la_j(\Om_0)<\la\}.
\end{equation}
and
\begin{equation}\label{nnormal}
n(\la,\Om_0)=
\begin{cases}
\underline{n}(\la,\Om_0)& \text{if }\la\not\in \sigma\big(H(\Om_0)\big)\\
\underline{n}(\la,\Om_0)+1& \text{if } \la\in\sigma\big(H(\Om_0)\big).
\end{cases}
\end{equation}
So we always have~:
\begin{equation}\label{nnn}
\underline{n} (\lambda, \Omega_0) \leq n (\lambda, \Omega_0) \leq
\overline{n}(\lambda,\Omega_0)\;.
\end{equation}
with equality when $\lambda$ is not an eigenvalue.
Note that $\overline{n}(\lambda,\Omega_0) - \underline{n} (\lambda, \Omega_0)$
    is the multiplicity of $\lambda$ when $\lambda$ is an eigenvalue of
$H(\Omega_0)$,
i.e. the dimension of the eigenspace associated to $\la$. We shall consider
    a family of $k$ open sets $\Omega_\ell$ \-($\ell=1,\dots,k$) 
contained in $\Omega_0$
    and the corresponding Dirichlet realizations $H(\Omega_\ell)$. For
each  $H(\Om_\ell)$ the corresponding eigenvalues counted with multiplicity are denoted by $(\lambda_k^\ell)_{k \in \mathbb N \setminus \{0\}}$ (with $\lambda_k^\ell \leq \lambda_{k+1}^\ell$). When counting the eigenvalues
less than some given $\lambda$ ,  we shall for simplicity write
\begin{equation}\label{notation}
n_\ell=n_\ell(\la)=n(\la,\Om_\ell)
\end{equation}
and analogously for the quantities with over-, respectively,
underbars.
\begin{theorem}\label{MAIN}~\\
Suppose $\Om_0$ is connected and that $\la\in\sigma\big(H(\Om_0)\big)$.
Suppose that the sets  $\Om_\ell$~($\ell= 1,\dots,k$) are pairwise disjoint
open subsets of $\Om_0$.
Then
\begin{equation}\label{mmm}
\sum_{\ell=1}^k \overline n_\ell\le n_0 +
\min_{\ell\ge 0}\big( \overline n_\ell-n_\ell\big)
\end{equation}
\end{theorem}
A direct weaker consequence of (\ref{mmm}) is the more standard
\begin{corollary}\label{Corw}~\\
Under the assumptions  of Theorem \ref{MAIN}, we have
\begin{equation}\label{mmw}
\sum_{\ell=1}^k \overline n_\ell\le \overline n_0\;.
\end{equation}
\end{corollary}
This corollary is actually present in the proofs of
    the asymptotics of the counting function (see for example
    the Dirichlet-Neumann bracketing in Lieb-Simon \cite{LiSi}).

\begin{remark}~\\
{\rm Inequality (\ref{mmm}) is also true if
$\la\not\in\sigma\big(H(\Om_0)\big)$.
The statement becomes more simply
$$
\sum_{\ell= 1}^k  \overline{n}_\ell \leq n_0 \;.
$$
and is proved essentially in the same way.}
\end{remark}
\begin{remark}~\\
{\rm The assumption that $\Om_0$ is connected is necessary. Indeed, suppose
$\Om_1$ and $\Om_2$ are connected and
assume that $\Om_0=\Om_1\cup\Om_2$ with $\Om_1\cap\Om_2=\emptyset$ and that
$\lambda=\la_1(\Om_1)=\la_1(\Om_2)$. Then $\la_1(\Om_0)=\la_2(\Om_0)$
and we deduce\break
$n(\lambda,\Om_0)=1$. In general we just have,  if we no longer assume the
connectedness of $\Om_0$,  Corollary \ref{Corw}.}
\end{remark}

Finally we  show that {\bf Courant's nodal Theorem} is an easy corollary of
Theorem \ref{MAIN}.

\begin{corollary}\label{Courant}{\bf :~Courant's nodal Theorem}\\
If $\Om_0$ is connected and if  $u$ is an eigenvector of $H(\Om_0)$
associated to  some eigenvalue $\la$, then
\begin{equation*}
\mu(u)\le n(\lambda,\Om_0)\;.
\end{equation*}
\end{corollary}
{\bf Proof}.\\
We now simply apply Theorem \ref{MAIN} by taking
$\Om_1,\dots, \Om_{\mu(u)}$  as the
nodal domains associated to $u$.
    We just have to use \eqref{u>0} for each $\Omega_\ell$,
$\ell=1,\dots, \mu(u)$,  which gives  $\overline n_\ell= n_\ell  =1$.
\finbox
\begin{remark}\label{litC}~\\
{\rm  Courant's nodal Theorem is one of the basic results in spectral
theory of Schr\"odinger type  operators. It is the natural generalization of
Sturm's oscillation theorem for second order ODE's. For recent investigations
see for instance \cite{A:1998} and \cite{D:2001}. }
\end{remark}
\section{Converse results.}\label{Section3}
In this section we consider some results converse to Theorem \ref{MAIN}.
\begin{theorem}\label{MAINc}~\\
Suppose that the $\Om_\ell$, $\ell\ge 1$, are pairwise disjoint
open  subsets of  $\Om_0$.
If $\la\in\sigma\big(H(\Om_0)\big)$ and
\begin{equation}\label{nnnc}
\sum_{\ell=1}^k \overline{ n_\ell} \geq n_0\;,
\end{equation}
then, for each $\Om_\ell$, $\la\in\sigma(H(\Om_\ell))$. If  $U_\ell(\la)$
denotes
the eigenspace  of $H(\Omega_\ell)$ associated to the  eigenvalue $\la$,
then there is an eigenfunction $u$ of $H(\Om_0)$ with eigenvalue
$\la$ such that
\begin{equation}\label{Phiell}
u=\sum_{\ell= 1}^k \varphi_\ell \mbox{ in } W_0^{1,2}(\Omega_0)\;,
\end{equation}
where each $\varphi_\ell$ belongs to $ U_\ell(\la)\setminus\{0\}$ and 
is identified with its extension by $0$ outside $\Omega_\ell$.
\end{theorem}
\begin{remark}~\\
One can naturally think that formula \eqref{Phiell} has immediate consequences
   on the family $\Omega_\ell$, which should have for example some covering
   property. The question is a bit more subtle
   because we do not want to assume a priori strong regularity properties
   for the boundary of the $\Omega_\ell$. We shall discuss this point
   in detail  in the last section.
\end{remark}

Another consequence of equalities in  Theorems \ref{MAIN} or
\ref{MAINc}
is given by the following results.
\begin{theorem}\label{subset}~\\
Suppose that, for some open subset $\Om_0$ in ${ \mathbb R } ^d$, some
$\lambda \in \sigma (H(\Omega_0))$  and some family of  pairwise
disjoint open
sets $\Om_\ell\subset \Om_0$, $0<\ell \leq k$, we have
\begin{equation}\label{mmmeq} 
\sum_{\ell= 1}^k \overline n_\ell =
n_0 +
\min_{\ell\ge 0}\big( \overline n_\ell-n_\ell\big)\;.
\end{equation}
Then,  for any subset   $L\subset \{ 1,2,\dots, k\}$,  such that
$\Om^*_L=\Inte\big(\cup_{\ell\in L}\overline{\Om_\ell}\big)\setminus \pa\Om_0$
    is connected, we have
\begin{equation}\label{m*}
\sum_{\ell\in L}\overline n_\ell=n(\la,\Om^*_L)+
\min\bigg(\min_{\ell\in L}\big(\overline n_\ell-n_\ell\big),
\: \overline n(\la,\Om^*_L)-n(\la,\Om^*_L)\bigg)\;.
\end{equation}
\end{theorem}
A simpler variant  is the following~:
\begin{theorem}\label{subseta}~\\
    Suppose \eqref{nnnc} holds and  that $\Om^*_L$ is defined as above. Then
    we have the inequality~:
\begin{equation}\label{n>n*}
\sum_{\ell\in L}\overline n_\ell\ge n(\la,\Om^*_L)\;.
\end{equation}
\end{theorem}

\paragraph{On the sharpness of Courant's nodal Theorem}~\\
It is well known that Courant's nodal Theorem is sharp only for
finitely many $k$'s \cite{Pleijel:1956}.

Let $\Omega_0$ be connected. We will say that an eigenvector $u$
    attached to an eigenvalue $\lambda$ of  $H(\Om_0)$ is {\bf Courant-sharp}
if $\mu(u) = n(\lambda,\Omega_0)$.
    Theorem \ref{subset} implies now~:
\begin{corollary}\label{submu}~\\
i) Let $u$ be a {\bf Courant-sharp} eigenvector of $H(\Omega_0)$ with
$\mu(u)=k$.
Let $\mathcal D^{(k)}=\{D_i\}_{i\in\{1,\dots,k\}} $
be the family  of the  nodal domains associated to $u$. Let $L$ be a subset of
$\{1,\cdots,k\}$  with $\# L =\ell$ and let
$\mathcal D_{L}$ be the  subfamily $\{D_i\}_{i\in L}$.
Let $\Om_L=\Inte (\overline{\cup_{i\in  L} D_i})\setminus\pa\Om_0$. Then
\begin{equation}\label{lk=ll}
\la_\ell(\Om_L)=\la_k
\end{equation}
where $\la_j(\Om_L)$ are the eigenvalues of $H(\Om_L)$.\\
ii) Moreover,
when $\Omega_L$ is connected, and if $\ell < k$, 
    $u\big |_{\Om_L}$ is {\bf Courant-sharp} and
$\lambda_\ell(\Omega_L)$ is simple.
\end{corollary}

\section{Basic tools}\label{Section4}
Let us first recall some  basic tools (see e.g. \cite{ReedSimon:1978})
    which were already  vital for
the proof  of  Courant's classical result.

\subsection{Variational characterization}
Let us first recall the variational characterization of eigenvalues.
\begin{proposition}\label{varchar}~\\
Let $\Omega$ be a bounded open set in $\mathbb R^d$ and $V$ real 
in $L^\infty(\Omega)$.
Suppose $\la  \in\sigma\big(H(\Om)\big)$ and let $
\mathcal U_{\pm}=\text{span }\langle u_1,\dots,
u_{k_\pm}\rangle $
where
\begin{equation}\label{kpm}
k_-=\underline{n}(\la,\Om)\text{ and } k_+=\overline{n}(\la,\Om)\;.
\end{equation}
Then
\begin{equation}\label{var}
\la =\inf_{\varphi\bot \mathcal U_-,\:\varphi\in W_0^{1,2}(\Om)}
\frac{\langle \varphi,\; H(\Om)\varphi\rangle}{\|\varphi\|^2}
\end{equation}
and
\begin{equation}\label{var+}
\la<\la_{\overline n(\la,\;\Om)+1}=\inf_{\varphi\bot\mathcal U_+,\:
\varphi\;\in W_0^{1,2}(\Om)}
\frac{\langle \varphi,\; H(\Om)\varphi\rangle}{\|\varphi\|^2}.
\end{equation}
If in \eqref{var} equality is achieved for some $\Phi\not\equiv 0$,  then
$\Phi$ is an eigenfunction in the eigenspace of $\la$.
\end{proposition}
Note that actually \eqref{var} and \eqref{var+} are the same statement.
We just state them separately for further reference. Note that
    we have not  assumed that $\Om$ is connected.

\subsection{Unique continuation}
   Next we restate a weak form of the unique continuation property:

\begin{theorem}\label{unique}~\\
Let $\Omega$ be an open set in $\mathbb R^d$ and $V$ real 
in $L^\infty_{loc}(\Omega)$.
Then any  distributional solution solution in $\Omega$ to $(-\Delta
+V)u=\la u$ which
vanishes on an open subset $\omega$ of $\Om$ is identically zero in
the connected component
    of $\Omega$ containing $\omega$.
\end{theorem}
There are stronger results of this type under weaker assumptions on the
potential, see \cite{JeKe}.
\subsection{ A consequence of Harnack's Inequality}
The standard Harnack Inequality (see e.g.\ Theorem 8.20 in \cite{GT:1983}),
together with the unique continuation theorem
    leads to the
\begin{theorem}\label{harnack}~\\
If $u$ is an eigenvector of $H(\Omega)$, then
for any $x$ in $N(u)\cap \Omega $ and any ball $B(x,r)$ ($r>0$), there
    exists $y_\pm \in  B(x,r)\cap \Omega $ such that $\pm u(y_\pm)>0$.
\end{theorem}

\section{Proof of the main theorems.}\label{Section5}
\subsection{Proof of Theorem \ref{MAIN}}
    Assume first for contradiction that
\begin{equation}\label{indirect}
\sum_{\ell\ge 1} \overline n_\ell>
n_0 +
\min_\ell \big( \overline n_\ell-
n_\ell \big)
\end{equation}
and recall that we assume  $\la\in\sigma\big(H(\Om_0)\big)$.
Pick some $\ell_0$ such that
$$
\overline n_{\ell_0}- n_{\ell_0}
=\min_\ell \big( \overline n_\ell- n_\ell \big) \;.
$$
{\bf Suppose first $\ell_0\geq 1$.}\\
We can rewrite \eqref{indirect} to obtain
\begin{equation}\label{indirect1}
\sum_{\ell\neq\ell_0,\;\ell\ge 1} \overline n_\ell +
n_{\ell_0}  > n_0\; .
\end{equation}

Let $\varphi^{\ell_0}_i,\:i=1,\dots,\underline n(\lambda,
\Omega_{\ell_0})$ denote the first
$\underline n_{\ell_0}$ eigenfunctions of $H(\Om_{\ell_0})$.
    The corresponding
eigenvalues are strictly smaller than $\la$. These functions and the remaining
$\sum_{\ell\neq\ell_0} \overline n_\ell $ eigenfunctions associated
to the other
$H(\Om_\ell)$ span a space of dimension at least $n_0$.
We can pick a linear combination $\Phi\not\equiv 0$ of these
functions which is orthogonal
to the $\underline n_0$  eigenfunctions of $H(\Om_0)$. By assumption
\begin{equation}\label{vPhi}
\frac{\langle\Phi, H(\Om_0)\Phi\rangle}{\|\Phi\|^2}\le \la,
\end{equation}
hence $\Phi$ must by the variational principle
be an eigenfunction and there must be equality in \eqref{vPhi}.

There are two possibilities: either  some
$\varphi^{\ell_0}_i, i<n_{\ell_0}$
    contributes to  the linear combination
which makes up $\Phi$ or
not. In the first case this means that the left hand side of \eqref{vPhi}
is strictly smaller than $\la$, contradicting the variational
characterization of $\la$. In the other case we obtain a contradiction
to unique continuation, since then $\Phi\equiv 0$ in $\Om_{\ell_0}$ and hence
vanishes identically in all of $\Om_0$.\\
{\bf Consider now the case when $\ell_0=0$.}\\
We have to  show that the assumption
\begin{equation}\label{n<nOm}
\sum_\ell \overline n_\ell>\overline n_0
\end{equation}
    leads to a contradiction. To this end it suffices to apply
\eqref{var+}. Indeed, we
can find a linear combination $\Phi$ of the eigenfunctions $\varphi^\ell_j$,
    $j\le \overline n_\ell$, 
corresponding to the different $H(\Om_\ell)$ such that
$0\not\equiv \Phi\bot \mathcal U_+$, but satisfies
\begin{equation*}
\frac{\langle \Phi,\;H(\Om_0)\;\Phi\rangle}{\|\Phi\|^2}\le 
\la= \la_{\overline n_0}\;,
\end{equation*}
and this contradicts \eqref{var+}. This proves \eqref{mmm}.

\subsection{ Proof of Theorem \ref{MAINc}}
The inequality \eqref{nnnc} implies that we can find  a non zero
$u \bot \mathcal U_-$ in
the span of the eigenfunctions $\varphi_j^\ell$, $j=1,\dots \overline n_\ell$,
of the different $H(\Om_\ell)$. Again by the variational characterization,
\eqref{var} and  \eqref{vPhi}  hold and hence $u$ must
be an eigenfunction.
\finbox
\subsection{Proof of Theorem \ref{subset}}
We assume (\ref{mmmeq}).
Without loss we might assume that we have labeled the $\Om_\ell$ such that
$L=\{1,\cdots,K\}$, with $K\leq k$.
Let $n_*=n(\la,\Om^*_L)$. We apply Theorem \ref{MAIN}
   to the family $\Omega_\ell$ ($\ell \in L$) and replace
$\Omega_0$  by $\Omega_L^*$ and obtain~:
\begin{equation}\label{contrineq}
\sum_{1\le \ell\le K}\overline n_\ell\leq n_*+
\min\big(\overline n_*-n_*,\min_{1\le \ell\le K}(\overline n_\ell-n_\ell)\big)\;.
\end{equation}
    We assume for contradiction that
\begin{equation}\label{contr}
\sum_{1\le \ell\le K}\overline n_\ell<n_*+
\min\big(\overline n_*-n_*,\min_{1\le \ell\le K}(\overline n_\ell-n_\ell)\big)\;.
\end{equation}
This implies
\begin{equation}\label{contr1}
\sum_{1\le \ell\le K}\overline n_\ell<\overline n_*\;,
\end{equation}
and
\begin{equation}\label{contr2}
\sum_{1\le \ell\le K}\overline n_\ell<n_* +
\min_{1\le \ell\le K}(\overline n_\ell-n_\ell)\;.
\end{equation}
Theorem \ref{MAIN}, applied to
    the family $\Omega^*_L, \Omega_\ell$ ($\ell >K)$,
implies that
\begin{equation}\label{n**1}
\overline n_*+\sum_{K<\ell\le k}\overline n_\ell \le n_0+\min\big(
\overline n_0-n_0,\,\min_{K<\ell\le k}(\overline n_\ell-n_\ell)\big)\;,
\end{equation}
and
\begin{equation}\label{n**2}
    n_*+\sum_{K<\ell\le k}\overline n_\ell \le n_0 \;.
\end{equation}
By adding (\ref{contr1}) and (\ref{n**1}), we get~:
\begin{equation} \label{cont3}
\sum_{1\leq \ell\leq k} \overline n_\ell < n_0 + \min\big(
\overline n_0-n_0,\,\min_{K<\ell\le k}(\overline n_\ell-n_\ell)\big)\;.
\end{equation}
By adding (\ref{contr2}) and (\ref{n**2}), we obtain
\begin{equation}\label{cont4}
\sum_{1\leq \ell\leq k} \overline n_\ell < n_0+ \min_{1\leq \ell\leq
K}(\overline n_\ell-n_\ell)\;.
\end{equation}
The combination of (\ref{cont3}) and (\ref{cont4}) is in contradiction with
(\ref{mmmeq}).

\subsection{Proof of Theorem \ref{subseta}}
For the case that \eqref{nnnc} holds \eqref{n>n*} can be shown similarly.
    \eqref{nnnc}
reads
\begin{equation*}
\sum_{1\le \ell\le k}\overline n_\ell\ge n_0\;.
\end{equation*}
We assume for contradiction that
\begin{equation}\label{contra}
\sum_{1\le \ell \le K}\overline n_\ell <n_*\;,
\end{equation}
where $n_*$ is defined as above. The addition of  (\ref{n**2})
    and (\ref{contra}) leads to a contradiction.
\finbox

\section{Illustrative examples}\label{Section6}
\subsection{Examples for a rectangle}
We illustrate Theorem \ref{MAIN} by the analysis of various examples
in rectangles.
Pick a rectangle   $\Om_0=(0,2\pi)\times (0,\pi)$ and take
$\Om_1=(0,\pi)\times (0,\pi)$ and consequently $\Om_2=(\pi,2\pi)\times(0,\pi)$.
The eigenvalues corresponding to $\Om_0$ for $-\Delta$  with Dirichlet boundary
conditions are given by
\begin{equation}\label{specrec}
\sigma\big(H(\Om_0)\big)=\bigg\{\la\in\mathbb R\:\bigg|\: \la=m^2/4+n^2,\;
(m,n)\in \mathbb Z^2,\: m,n>0\bigg\}\;,
\end{equation}
while those for $\Om_1$ and hence for $\Om_2$ which can be obtained by
a translation of $\Om_1$.
are given by
\begin{equation}\label{specsquare}
\sigma\big(H(\Om_1)\big)=\sigma\big(H(\Om_2)\big)=
\bigg\{\la\in\mathbb R\:\bigg|\: \la=m^2+n^2,\;(m,n)\in \mathbb Z^2,\:
m,n>0\big\}.
\end{equation}
Denote the eigenvalues associated to $\Om_0$ by $\{\la_i\}$ and those to
$\Om_1$ by $\{\nu_i\}$. We easily check that
$\la_5=\la_6=\nu_2=\nu_3=5$, $\la_{11}=\la_{12}=\nu_5=\nu_6=10$ so that
for these cases Theorem \ref{MAIN} is sharp.

One could  ask whether  there are arbitrarily  high eigenvalues
cases for which we have equality in \eqref{mmm}.
This is not the case,
as can be seen from the following standard number theoretical considerations.
We have (see \cite{VdK} and for more recent contributions
    \cite{Ran} and \cite{Ber})  the following asymptotic estimate
for the number of lattice points in an ellipse.
Let $a,b>0$, then
\begin{equation}\label{lattice}
A(\la):=\#\bigg\{(m,n)\in\mathbb Z^2\:\bigg|\:am^2+bn^2\le \la\bigg\}
\end{equation}
has the following asymptotics as $\la$ tends to infinity:
\begin{equation}\label{asympt}
A(\la)=\frac{\pi}{\sqrt{ab}}\la+\Og(\la^{1/3}).
\end{equation}
We have not to consider $A(\lambda)$ but rather
\begin{equation}\label{A+}
A^+=\#\bigg\{(m,n)\subset\mathbb Z^2, m,n>0\:\bigg|\:am^2+b n^2\le \la\bigg\}\;.
\end{equation}
Hence we get
\begin{equation}\label{A+A}
\begin{split}
A(\la)=4A^+(\la)+2\#\bigg\{m\in \mathbb N,\; m>0\;\bigg|\;m
\le \big[(\la/a)^{1/2}\big]\bigg\}\\
+2\#\bigg\{n\in \mathbb N,\;n>0\;\bigg|\;n\le\big[(\la/b)^{1/2}\big]\bigg\}+1\;.
\end{split}
\end{equation}
If we apply this to $A^+$ with $a=1/4, b=1$ (in this case denoted by
$A_0^+$) and to $A^+$ with $a=1, b =1$ (in this case denoted by
$A_1^+$),
we get asymptotically
\begin{equation}\label{A-A}
A_0^+(\la)-2A_1^+(\la)= \frac{1}{2}\sqrt\la+o\,(\sqrt \la)\;.
\end{equation}
Note that
$$\overline n_i(\la)= A_i^+(\la),\;i=0,1\;.$$
In order to  control $n_i(\lambda)$, we observe that, for any $\epsilon >0$~:
$$\overline n_i(\la - \epsilon)\leq n_i(\lambda)\leq \overline n_i(\la)\;.$$
This implies
\begin{equation}\label{logbis}
\overline n_i(\la)-n_i(\la)=\Og (\la^\frac 13)\;.
\end{equation}
The asymptotic formula \eqref{asympt} implies
\begin{equation}\label{log}
\overline n_i(\la)-n_i(\la)=o(\sqrt\la)\;,
\end{equation}
and this shows that for large $\la $ \eqref{mmm} is never sharp.

\subsection{About Corollary \ref{submu}.}
One can ask
whether there is a converse to Corollary \ref{submu}  in the following
sense.  Suppose we have an eigenfunction $u$ with $k$ nodal domains and
eigenvalue $\la$. For each pair of neighboring nodal
domains of $u$, say, $D_i$ and $D_j$, let
$\Om_{i,j}=\Inte\;(\overline{D_i\cup D_j})$ and suppose that
$\la=\la_2(\Om_{i,j})$.
Does this imply  that $\la=\la_k$? The answer to the question is negative,
as the following easy example shows~:\\
Consider the rectangle $Q=(0,a)\times(0,1)\subset \mathbb R^2$ and consider
$H_0(Q)$. We can work out the eigenvalues explicitly as
\begin{equation}\label{Qmn}
\{\pi^2(\frac{m^2}{a^2}+n^2)\},\: \text{ for } m,n\in \mathbb N\setminus 0,
\end{equation}
with corresponding eigenvectors
$(x,y)\mapsto \sin (\pi m \frac xa )(\sin \pi n y)  $.
If
\begin{equation}\label{<a>}
a^2\in \big(\frac{9}{4},\:\frac{8}{3}\big)\;,
\end{equation}
then
\begin{equation*}
\la_3(Q)=\pi^2(\frac{1}{a^2}+4)<\la_4(Q)=\pi^2(\frac{9}{a^2}+1)\;,
\end{equation*}
and the zeroset of $u_4$ is given by $\{(x,y)\in Q\:|\: x=a/3,\: x=2a/3\}$.
For $u_4$ we have  $\Om_{1,2}=Q\cap \{0<x<2a/3\}$. If $2a/3>1$ (which is
    the case under  assumption \eqref{<a>}),  then $\la_2(\Om_{1,2})
=\la_4(Q)$. We have consequently an example with $k=3$.
 
\section{Converse theorems in the case of regular open sets}\label{ancona}
\subsection{Preliminary discussion about regularity}\label{preldisc}
Before we present what we think would be the right notion of regular
   open set adapted to our problem, let us discuss briefly other
possible notions.\\
    As a consequence of Theorem \ref{MAINc} and  using \eqref{u>0},
we get that  each nodal
domain $D_{k\ell}$  of $\varphi_\ell$
    is included in  a nodal domain $D_{j0}$ of $u$. Using a result of
Gesztesy and Zhao  (\cite{GZ:1994}, Theorem 1), this implies also that the
capacity (see next subsection)
    of $D_{j0} \setminus D_{k\ell}$
(hence  the measure)
    is $0$.\\
    At the ``regular'' points of the boundary of  $\Omega_\ell$  one
can get additional information.\\
    Let us say that,
if $\Omega$ is  an open set and if
$\pa \Omega = \overline{\Omega} \setminus
    \Omega$,  that a point of $\pa \Omega$ is $C^{1,\alpha}$-regular if
    there exists a neighborhood $\mathcal V(x)$ of $\pa \Omega$
    such that $\pa \Omega \cap \mathcal V(x)$ is a $C^{1,\alpha}$- hypersurface.
    We denote by $\pa \Omega ^{c1reg}$ the subset of the regular points.\\
 
Then  we get, from \eqref{Phiell}, using the property of the
restriction  map from $W^{1,2}(\Omega_\ell)$ into $W^{\frac 12,2}
(\partial  \Omega_\ell ^{c1reg})$, that,
under the assumptions of Theorem \ref{MAINc}, we have the inclusion
\begin{equation}\label{Nu}
\bigcup_{\ell\ge 1} (\pa\Om_\ell^{c1reg})  \subset
N(u) \cup \pa \Omega_0\;.
\end{equation}

    One could say that an open set $D$ is topologically regular
    if the subset of regular points is dense in $\pa D$.\\
In this case, it is easy to see that if
some function $f$ belongs to 
$ C^0(\overline{D}\cap \Omega_0)\cap W^{1,2}_0(D)$,  then
    $f$ vanishes on the boundary (we first get it at the regular points
    and then conclude by continuity).\\
More generally, one could ask under which weakest condition on a point
    $x$ of $\pa D\cap\Omega_0$ any function $f$
in   $C^0(\overline{D}\cap \Omega_0)\cap W^{1,2}_0(D)$ satisfies
$f(x)=0$. This is what we will discuss in the next subsections.

\subsection{Capacity}\label{ss.cap}
There are various equivalent definitions of polar sets and capacity
(see e.g.\\ \cite{Den}, \cite{Fu}, \cite {He}, \cite {LiLo}). If $U$ is a
bounded open subset
of ${ \mathbb R } ^d$, we denote by $\Vert .\Vert_{W^{1,2}_0(U)}$ the Hilbert
norm on $W_0^{1,2}(U)$~:
$$u\mapsto \Vert u\Vert_{W^{1,2}_0(U)} :=(\int_U \vert \nabla u\vert
^2\, dx)^{\frac 1 2}\;.$$
   The
capacity in $U$ of $A\subset U$  can be defined\footnote {For $d\geq 
3$ the restriction
that $U$ is bounded can be removed and one may take $U={ \mathbb R }
^d$.} as   $$
\begin{array}{ll}
\text {
Cap}_U(A) \;
&:=\;\inf \{ \Vert s\Vert
^2_{W^{1,2}_0(U)}\,;\,s\in W^{1,2}_0(U)\quad\quad \\
&\quad\quad \quad\quad\quad  \text { and } s\geq 1  \text{ 
a.e.\ in  some
neighborhood of } A\,\} \;. 
\end{array}$$ It is easily checked
that if $K$ is  compact and $K\subset U\cap V$, where
$V$ is also  open and bounded in ${ \mathbb R } ^d$, there is a
$c=c(K,U,V)$ such that $\text {Cap}_U(A)\leq c\; \text {Cap}_V(A)$
for $A\subset K$.
So $\text {Cap}_U(A)=0$ for some bounded open $U\supset A$ iff for
each $a\in A$ there exists
an $r>0$ and a bounded region $V$ such that $V\supset B(a,r)$ and
$\text {Cap}_{V}(B(a,r)\cap A)=0$.
In this case we may simply write $\text {Cap}(A)=0$ without referring to $U$.

\subsection{Converse theorem}\label{convthms}
We are now able to formulate our definition of regular point.
\begin{definition}~\\
Let $D$ be an open set in ${ \mathbb R } ^d$.
We shall say that a point $x\in \pa D$
    is (capacity)-regular (for $D$) if, for any $r >0$,  the capacity
of $B(x,r)\cap \complement D$
    is strictly positive.
\end{definition}

\begin{theorem}\label{capreg}~\\
Under the assumptions of Theorem \ref{MAINc}, any
point $x\in
\partial \Omega_\ell \cap \Omega_0$ which is (capacity)-regular with
respect to $\Omega _\ell$ (for some $\ell$) is in the nodal set of
$u$.
\end{theorem}
This theorem admits the
\begin{corollary}~\\
Under the assumptions of  Theorem \ref{MAINc} and if, for all $\ell$,
every point in
    $(\partial \Omega_\ell )\cap \Omega_0$ is (capacity)-regular for
$\Omega _\ell$, then the family of the nodal domains of $u$
    coincides with the union over $\ell$ of the family of the nodal 
domains of the
$\varphi_\ell$, where $u$ and $\varphi_\ell$ are introduced in \eqref{Phiell}.
\end{corollary}
{\bf Proof of corollary}\\
It is clear that any nodal domain of $\varphi_\ell$ is contained in 
contained in a unique  nodal domain of $u$.\\
Conversely, let $D$ be a nodal domain of $u$ and let $\ell \in 
\{1,\dots,k\}$. Then, by combining
  the assumption on $\pa \Omega_\ell$, Proposition \ref{capzero}
  and \eqref{Phiell}, we obtain the property~:
$$
\pa \Omega_\ell \cap D =\emptyset\;.
$$
Now,  $D$ being connected,  either $\Omega_\ell \cap D=\emptyset$ or 
$D\subset \Omega_\ell$. Moreover the second case should occur for at 
least one $\ell$,
say $\ell=\ell_0$. Coming back to the definition of a nodal set and 
\eqref{Phiell}, we observe that  $D$
    is necessarily contained in a nodal domain $D_j^{\ell_0}$ of 
$\varphi_{\ell_0}$.\\
Combining the two parts of the proof gives that any nodal set of $u$ 
is a nodal set of $\varphi_\ell$ and vice-versa.
 
\subsection{Proof of Theorem \ref{capreg}}
According to the discussion of Subsection \ref{preldisc}, the proof will be
    a consequence of the following proposition.
\begin{proposition}\label{capzero}~\\
Let us consider two open subsets $D$ and $\Omega$  of $\rz^d$ such
that $D \subset \Omega $ and
    a point $x_0$  in $\partial D \cap \Omega$.
Assume that, for some given $r_0 >0$ such that $B(x_0,r_0)\subset \Omega$,
    there exists $u\in W^{1,2}_0(D)$ and $v \in  C^0 ( B(x_0,r_0))$
    such that~:
$$
    u_{\vert D\cap  B(x_0,r_0)} = v_{\vert  D\cap  B(x_0,r_0)}\; \text{
a.e. in } D\cap  B(x_0,r_0)\;.
$$
Then if $ v (x_0) \neq 0$, there exists
    a ball $B(x_0,r_1)$ ($r_1 >0$),  such that $B(x_0,r_1)\setminus   D$
is polar, that is of
capacity $0$. \end{proposition}

\begin{remark}~\\ Using some standard  potential theoretic arguments, 
Proposition \ref{capzero} can be deduced  from Th\'eor\`eme 5.1 in
     \cite{DenL} which characterizes those $u\in W^{1,2}(\Omega )$ 
that belong to  $W_0^{1,2}(\Omega )$. The proof below should be more 
elementary in character.\\
\end{remark}

\begin{remark}~\\
    Given a region $D\subset { \mathbb R } ^n$ and a ball $B=B(x,r)$, 
$x\in \partial D$, the difference set $B(x,r)\setminus   D$
is polar if and only if  $B(x,r)\cap  \partial D$
    is polar. This follows from the fact that  a polar
subset of $B=B(x,r)$ does not disconnect $B$ \cite{Bre}.\\

\end{remark}
\begin{remark}~\\
If $D$ is a nodal domain of an eigenfunction $u$ of $H(\Omega)$, then
any point of $\partial D
    \cap \Omega$ is capacity-regular for $D$. This is an immediate
consequence of Theorem~\ref{harnack} (it also follows from the 
preceding remark).
Indeed, if $x$ is  in $\pa D \cap \Omega$, then for
any $r>0$, one can find a ball $B(y,r')$ in $\complement D\cap
B(x,r)$.
\end{remark}

To prove Proposition \ref {capzero} we require some well-known facts stated in
the next three lemmas.
\begin{lemma}\label {lem1}~\\
Let  $U$ be a bounded convex domain  in ${ \mathbb R } ^d$ and let
$B(a,\rho )$, $\rho >0$ be a ball such that $\overline
B(a,\rho )\subset U$. There exists a positive constant  $c=c(a,\rho
,U)$ such that, for every $f\in W^{1,2}(U)$ vanishing a.e.\ in
$B(a,\rho )$,
$$\Vert f\Vert _{W^{1,2}(U)} \leq c\, \Vert \nabla f\Vert _{L^2(U)}\;.$$
\end{lemma}

\paragraph{\bf Proof.}~\\
   We assume as we may that $a=0$ and let $U'=U\setminus B(0,\rho
)$. Fix $R$ so large that $U\subset B(0,R)$. By approximating $f$ by 
smooth functions (e.g.\ regularize
$f((1-\delta
){\mathbf .})$ for $\delta >0$ and small to get $f_1\in C^\infty
(\overline U)$), we may restrict to functions $f\in C^\infty
(U)$ vanishing in $B(0,\rho )$. Then, since $$ \vert
f(x)\vert ^2=\vert   \int
_0^1\,\nabla f(sx).x\, ds\vert  ^2\leq R ^2  \int _{\frac{\rho}{ \vert
x\vert }}^1\,\vert \nabla f(sx)\vert ^2\, ds \text{ for }x\in U'\;,$$
 we have
\begin{equation}\label{NNN}
\begin{split}
\int_{U'} \vert f(x)\vert ^2\, dx &\leq R^2\,\iint _{x\in 
U',{\frac{\rho}{   \vert x\vert} }\leq s\leq 1} \vert \nabla 
f(sx)\vert ^2\, dx\,
ds\\
&\leq R^2\iint_{z\in sU',\,{\rho \leq  \vert z\vert },\, s\leq 1}
\vert \nabla f(z)\vert ^2\, dz\, {\frac{ds}{s}}\\
& \leq \frac{R^3}{ \rho } \; \int _{U'}\, \vert \nabla f(x) \vert ^2\; dx, \\
\end{split}
\end{equation}
and the lemma follows.
\begin{lemma}\label{lem2}~\\
   Let $U$ be a domain in ${ \mathbb R } ^d$. For every $f\in W^{1,2}(U)$
the function $g=f_+ $  is also in $W^{1,2}(U)$, with $\Vert g\Vert
_{W^{1,2}(U)}\leq \Vert f\Vert _{W^{1,2}(U)}$. Moreover the map $f\mapsto g$
from $W^{1,2}(U)$ into itself is continuous (in the norm topology).
\end{lemma}

\begin{remark}~\\
Since $\inf \{
f_n,1\}=1-(1-f_n)_+$, it follows from the lemma that $\inf \{ f_n,1\}
\to \inf \{ f,1\}$ in $W^{1,2}(U)$ whenever $f_n\to f$ in
$W^{1,2}(U)$.
\end{remark}

\paragraph{\bf Proof.}~\\
For the first two facts we refer to \cite {KiSt} or \cite {LiLo},
where it is moreover shown that the weak partial derivatives
$\partial _{j}{f_+}$
and $\partial _{j}f$ satisfy
$$\partial _{j}{f_+}=1_{\{f>0\}}\,\partial _{j}f=1_{\{f\geq 
0\}}\,\partial _{j}f \text { \ \ a.e. in \ \ } U.$$
Therefore, for any $\delta >0$, we have~:
\begin{equation}\label{BBB}
\begin{split}
\Vert \nabla [f_n]_+-\nabla f_+\Vert _{L^2}&=\Vert 1_{\{ f_n>0\} }
\nabla f_n-1_{\{ f>0\} }\nabla f\Vert _{L^2}\\
   &\leq \Vert 1_{\{ f_n>0\} } (\nabla f_n-\nabla
f)\Vert _{L^2}+\Vert (1_{\{ f>0\} }-1_{\{ f_n>0\} })\nabla f\Vert _{L^2} \\
& \leq \Vert \nabla f_n-\nabla
f\Vert _{L^2}+\Vert (1_{\{ f>0;f_n\leq 0\} }+1_{\{ f\leq 0;f_n> 0\} 
})\nabla f\Vert _{L^2}\\
   &\leq \Vert \nabla f_n-\nabla
f\Vert _{L^2}+\Vert 1_{\{ 0\leq \vert f\vert \leq  \delta\} } \nabla 
f\Vert _{L^2}+2\Vert 1_{\{ \vert f_n-f\vert \geq \delta \} }\nabla 
f\Vert
_{L^2}.
\end{split}
\end{equation}

Given  $\varepsilon >0$, fix $\delta >0$ so that $\Vert 1_{\{ 0\leq 
\vert f\vert \leq
\delta\} } \nabla f\Vert _{L^2}\leq \varepsilon $ (recall that 
$\nabla f=0$ a.e.\ in $\{ f=0\} $). Since
$\nabla f\in
L^2(U )$ and  $\Vert 1_{\{ \vert f-f_n\vert \geq \delta \} }\Vert 
_{L^1}\leq \frac{\Vert
f_n-f\Vert _{L^2}^2}{ \delta ^2}$,  it follows that $\displaystyle \lim_{n\to
\infty } \Vert
(1_{\{ \vert f-f_n\vert \geq \delta \} })\nabla f\Vert _{L^2}=0$. Therefore
$\displaystyle \limsup_{n\to \infty }  \Vert \nabla [f_n]_+-\nabla
f_+\Vert _{L^2}\leq \varepsilon $, 
which proves  that $[f_n]_+\to f_+$ in $W^{1,2}(U)$,  
if $f_n\to f$ in $W^{1,2}(U)$.

\begin{lemma}\label{lem3}~\\ Let $\omega$ be open in ${ \mathbb R }^d$ and
let $\{ f_n\} $ be a sequence of continuous functions in $\omega$ such that
$f_n\in W^{1,2}(\omega)$ for each $n\geq 1$ and $\displaystyle \lim_{n\to
\infty }\Vert f_n\Vert _{W^{1,2}(\omega)}= 0$.\\
Then the set $\displaystyle F=\{ x\in \omega\,;\, \liminf_{n\to \infty }
\vert f_n(x)\vert >0\, \}
$ is polar.

\end{lemma}

{\bf Proof.}\\ It suffices to show that $\text {cap}_\omega(F\cap K)=0$ 
for any compact subset $K$ of $\omega $. Let $\varphi \in C_0^\infty 
({ \mathbb R }^d
)$ be such that
$0\leq \varphi \leq 1$ in ${ \mathbb R }^d$, $\varphi =1$ in $K$ and
$\text {supp}(\varphi )\subset \omega$. Then $g_n=f_n\varphi \to 0$ in
$W_0^{1,2}(\omega)$ and $g_n=f_n$ in
$K$.

   Set
$F_\nu =\{ x\in \omega\,;\, \vert g_n(x)\vert \geq 2^{-\nu }$ for all
$n\geq \nu \} $. By the definition of the capacity, we have  $\text {Cap}_\omega (F_\nu
)\leq 2^{2\nu} \Vert \nabla
g_n\Vert _{L^2}^2$ for all $n\geq \nu $ and $\text {cap} (F_\nu )=0$.
Therefore $\text {cap}_\omega (\bigcup _{\nu \geq 1}F_\nu )=0$ and $\text
{cap}_\omega (F\bigcap K )=0$, 
since $F\bigcap K\subset \bigcup _{\nu \geq 1} F_\nu $.

\begin{proposition}\label{prop1}~\\ Let $U$ be a non empty open subset of the
ball $B=B(a,r)$ in ${ \mathbb R } ^d$. Suppose there exist a continuous 
function $f$ in $U$
and a  sequence $\{ f_n\}
$ of continuous functions in  $B$ such that \\
(i) $ f\geq  1$ in $U$ and $f\in
W^{1,2}(U)$, \\
(ii)  $f_n=0$ in a neighborhood of $B\setminus U$ and $f_n\in
W^{1,2}(U)$ for each $n\geq 1$,\\
   (iii) $\displaystyle \lim_{n\to \infty } \Vert
f-f_n\Vert _{W^{1,2}(U)}=0$.\\
Then the set $F:=B\setminus U$ is polar.
\end{proposition}

\paragraph{Proof of Proposition \ref{prop1}}~\\
Replacing $f$ by $\inf \{ f,1\} $ and $f_n$ by
$\inf \{ f_n,1\} $, we see\footnote {The  weak
convergence $\inf \{ f_n,1\}\buildrel w \over \to \inf \{
f,1\} $ suffices here. It allows the approximation of $1=\inf \{ f,1\}
$ in the norm topology in $W^{1,2}(U)$ by finite convex combination
of  the $\inf \{ f_n,1\}$. So we are again left with the case when
$f=1$ in $U$.} from Lemma \ref{lem2}  that we may
assume that $f=1$ in $U$.

   So $\displaystyle \lim _{n\to
\infty }\Vert \nabla f_n\Vert _{L^2(U)}=0$ and $\displaystyle
\lim_{n\to \infty } \Vert 1-f_n\Vert _{L^2(U)}=0$.

   Fix a ball $\overline B(z_0,2\rho )\subset U$, $\rho >0$, and a
cut-off function $\alpha \in C^\infty ({ \mathbb R } ^d)$ such that
$\alpha =1$ in  $B(z_0,\rho )$,
$\alpha =0$ in ${ \mathbb R } ^d\setminus B(z_0,2\rho )$. Set 
$g=1-\alpha $, $g_n=(1-\alpha )f_n$.

Then $g$, $g_n$ belong to $W^{1,2}(B)$, $\nabla g=\nabla g_n=0$ a.e.\
in $F$ and
   $$\displaystyle \lim_{n\to \infty } \Vert \nabla (g-g_n)\Vert
_{L^2(B)}=\lim_{n\to \infty }
\Vert \nabla (g-g_n)\Vert _{L^2(U)} =0.$$

So, by Lemma \ref {lem1}, $\displaystyle \lim_{n\to \infty } \Vert
g-g_n\Vert _{W^{1,2}(B)}=0$. But  $g-g_n\geq 1$ in $ F$ and it
follows from Lemma \ref {lem3}
that $F$ is polar.

\paragraph{Proof of Proposition \ref{capzero}}~\\ Choose $r_1>0$ so
small that $v\geq c_0:=\frac 12 v(x_0)$ in $B(x_0,r_1)$. Since
$u\in W^{1,2}_0(D)$,  there
is a sequence $\{ u_n\} $ in $C_0^\infty ({ \mathbb R } ^d)$ such that
$\text {supp} (u_n)\subset D$ and $u_n\to u$ in $W^{1,2}({ \mathbb R }
^d)$. Applying
Proposition \ref{prop1}  to the ball $B(x_0,r_1)$ and the functions
$f=c_0^{-1}u_{\vert B(x_0,r_1)}$, $f_n=c_0^{-1}{u_n}_{\vert
B(x_0,r_1)}$, we see that $B(x_0,r_1)\setminus   D$
is polar.

{\bf Acknowledgement.}\\  It is a pleasure to thank
    Fritz Gesztesy and
Vladimir Maz'ya for useful discussions. The two last authors were partly
   supported by the European Science Foundation Programme {\it Spectral Theory
   and Partial Differential Equations} (SPECT) and the EU IHP network {\it
   Postdoctoral Training Program in Mathematical Analysis of Large
Quantum Systems} HPRN-CT-2002-00277.

\footnotesize
\bibliographystyle{plain}

\begin{thebibliography}{1}
\bibitem{A:1998}
G.~Alessandrini.
\newblock On Courant's nodal domains theorem.
\newblock{\em Forum Math.}, 10, p.~521-532 (1998).

\bibitem{Ber} P.~B\'erard.
\newblock Lattice points in some domains.
\newblock Comm. Partial Differential Equations, Vol 3 (4) p.~335-348 (1978).

\bibitem{Bre} M.~Brelot.
\newblock {\it \'El\'ements de la th\'eorie
classique du Potentiel.}
\newblock Les cours de Sorbonne. \newblock Centre de Documentation
Universitaire, Paris 1965.

\bibitem{D:2001}
E.B. Davies, G. Gladwell, J. Leydold, and P. Stadler.
\newblock Discrete nodal domain theorems.
\newblock {\em Linear Algebra Appl.} 336,  p.~51-60 (2001).

\bibitem{Den} J.~Deny.
\newblock CIME Potential theory.  Stresa, 2-10 July 1970, p.~121-201.
\newblock Edizioni Cremonese Roma 1970.

\bibitem{DenL} J.~Deny, J.~L.~Lions.
\newblock Les espaces du type de Beppo Levi.
\newblock {\em Ann. Inst. Fourier, Grenoble} 5, p.~305-370 (1955).

\bibitem{Fu} M.~Fukushima, Y.~Oshima, and  M.~Takeda.
\newblock Dirichlet forms and symmetric Markov processes.
\newblock {\em de Gruyter Studies in Mathematics}, de Gruyter \& Co.,
Berlin 1994.

\bibitem{GZ:1994}
F.~Gesztesy, Z.~Zhao.
\newblock Domain perturbations, Brownian motion, and groundstate
of Dirichlet Schr\"odinger operators.
\newblock{\em Math. Z.} 215, p.~143-150 (1994).

\bibitem{GT:1983}
D.~Gilbarg, N.S. Trudinger.
\newblock Elliptic Partial Differential Equations of Second Order.
\newblock{\em Springer} 1983.

\bibitem{He}
R. Helms.
\newblock Introduction to potential theory
\newblock{\em Wiley-Interscience, New York-London-Sydney,}  1969.

\bibitem{JeKe} D.~Jerison, C.E.~Kenig.
\newblock Unique continuation and absence of positive eigenvalues for 
Schr\"odinger operators (with an appendix by E.M.~Stein).
\newblock Annals of Math. (2) 121, n$^o3$, p. 463-494  (1985).

\bibitem{KiSt} D.~Kinderlehrer, G.~Stampacchia.
\newblock {\it An introduction to variational inequalities and their
applications.}
\newblock Pure and Applied Mathematics 88. Academic Press.  New
York-London, 1980.


\bibitem{LiLo} E.~Lieb, M.~Loss.
\newblock  {\it Analysis.} \newblock Graduate Studies in Mathematics, 14.
\newblock American Mathematical Society, Providence, RI, 1997.

\bibitem{LiSi} E.~Lieb, B.~Simon.
\newblock The Thomas-Fermi theory of atoms, molecules and solids,
\newblock {\it Adv. in Math.} 23, p.~22-116 (1977).

 
\bibitem{Pleijel:1956} A.~Pleijel.
\newblock Remarks on Courant's nodal theorem.
\newblock{\em Comm. Pure. Appl. Math.} 9, p.~543-550 (1956).

\bibitem{Ran} B.~Randol.
\newblock A lattice point problem.
\newblock Trens. AMS 121, p.~257-268 (1966).

\bibitem{ReedSimon:1978}
M.~Reed, B.~Simon.
\newblock Methods of modern mathematical physics IV: Analysis of operators.
\newblock {\em Academic Press}, 1978.

\bibitem{VdK}J.~Van der Corput.
\newblock Zahlentheoretische Absch\"atzungen mit Anwendungen auf
    Gitterpunkteprobleme.
\newblock Math. Zeitschrift 17, p.~250-259 (1923).

\end{thebibliography}

\scshape
A. Ancona: D\'epartement de Math\'ematiques, Bat. 425,
Universit\'e Paris-Sud, 91 405 Orsay Cedex, France.

email: Alano.Ancona@math.u-psud.fr

\scshape
B. Helffer: D\'epartement de Math\'ematiques, Bat. 425,
Universit\'e Paris-Sud, 91 405 Orsay Cedex, France.

email: Bernard.Helffer@math.u-psud.fr

\scshape
T. Hoffmann-Ostenhof: Institut f\"ur Theoretische Chemie, Universit\"at
Wien, W\"ahringer Strasse 17, A-1090 Wien, Austria and International Erwin
Schr\"odinger Institute for Mathematical Physics, Boltzmanngasse 9, A-1090
Wien, Austria.

email:thoffman@esi.ac.at

\end{document}